\documentclass[a4paper, 11pt]{article}
\usepackage[english]{babel}
\usepackage{latexsym}
\usepackage{amsfonts}
\usepackage{amssymb}
\usepackage{amscd}
\usepackage{amsmath}
\usepackage{graphics}
\usepackage{eepic}
\usepackage{epsfig}
\usepackage{theorem}

\newtheorem{theorem}{Theorem}[section]

\newtheorem{lemma}[theorem]{Lemma}

\newtheorem{remark}[theorem]{Remark}

\newcommand{\C}{{\mathbb C}}
\newcommand{\Z}{{\mathbb Z}}
\newcommand{\fie}{\varphi}
\newcommand{\I}{\mathcal{I}}
\newcommand{\<}{\backslash}

\numberwithin{equation}{section}

\begin{document}
\setlength{\unitlength}{1mm}

\title{\textbf{\large{THE MONODROMY CONJECTURE FOR ZETA FUNCTIONS ASSOCIATED TO IDEALS IN DIMENSION TWO}}}
\author{Lise Van Proeyen and Willem Veys
\footnote{K.U.Leuven, Departement Wiskunde, Celestijnenlaan 200B,
B-3001 Leuven, Belgium, email: Lise.VanProeyen@wis.kuleuven.be,
Wim.Veys@wis.kuleuven.be. The research was partially supported by
the Fund of Scientific Research - Flanders (G.0318.06).}}
\date{}

\maketitle

\begin{abstract}
The monodromy conjecture states that every pole of the topological
(or related) zeta function induces an eigenvalue of monodromy.
This conjecture has already been studied a lot. However in full
generality it is proven only for zeta functions associated to
polynomials in two variables.
\\
\indent In this article we work with zeta functions associated to
an ideal. First we work in arbitrary dimension and obtain a
formula (like the one of A'Campo) to compute the `Verdier
monodromy' eigenvalues associated to an ideal. Afterwards we prove
a generalized monodromy conjecture for arbitrary ideals in two
variables.
\end{abstract}

 \begin{center} \footnotesize{\emph{2000 Mathematics Subject Classification.} 14E15,
32S40, 14H20.}
\end{center}

\section{Introduction}

Classically the invariants called topological, motivic and
$p$-adic Igusa zeta function are associated to \emph{one}
polynomial $f$ over $\C,$ over an arbitrary field of
characteristic zero and over a $p$-adic field, respectively. There
are fascinating conjectures relating their poles with the roots of
the Bernstein-Sato polynomial (also called $b$-function) of $f$
and with the eigenvalues of the local Milnor monodromy of $f$ in
points of $\{f=0\},$ up to now proven in full generality only for
polynomials in two variables.

One associates similarly in a natural way all these functions to
several polynomials or to an ideal. Not so obvious is the notion
of Bernstein-Sato polynomials associated to several polynomials or
an ideal. There is a construction of Sabbah \cite{Sabbah} and more
recently also of Budur, Musta\c{t}\v{a} and Saito \cite{Budur -
Mustata - Saito}. Concerning monodromy, the classical construction
of local Milnor fibre does not generalize to arbitrary maps $(f_1,
\ldots , f_r).$ However, there is a notion of `Verdier monodromy'
in this general context \cite{Verdier artikel}.

In this paper we prove a relation between the poles of these zeta
functions associated to an arbitrary ideal in two variables and
the `Verdier monodromy eigenvalues' of this ideal, generalizing
the result for one polynomial.

\vspace{\baselineskip}

We now provide more details, focussing on the topological zeta
function. Let $\I = (f_1, \ldots , f_r)$ be a nontrivial ideal in
$\C[x_1, \ldots , x_n]$ and $Y = \mathbf{V}(\I)$ its associated
subscheme of $\mathbb{A}^n_{\C}.$ We assume that $Y$ contains the
origin.

We first fix notation to define the topological zeta function.
Take a principalization $\psi : \tilde{X} \to \mathbb{A}^n$ of
$\I.$ By this we mean that $\psi$ is a proper birational map from
a nonsingular variety $\tilde{X}$ such that the total transform
$\psi^*\I$ is a principal ideal with support a simple normal
crossings divisor, and moreover that the exceptional locus of
$\psi$ is contained in the support of $\psi^*\I$. If $\mathcal{I}$
has components of codimension one, we can write this total
transform as a product of two (principal) ideals: the support of
the first one is the exceptional locus, where the support of the
second one is formed by the irreducible components of the total
transform that are not contained in the exceptional locus. This
second ideal is the `weak transform' of $\mathcal{I}.$

Note that we use the word `principalization', where other authors
may also use log-principalization, log-resolution or
monomialization.

Let $\tilde{E} = \sum_{i \in J}N_iE_i$ denote the divisor of
$\psi^*\I$, i.e.\ its irreducible components are the $E_i,$ $ i
\in J,$ occurring with multiplicity $N_i.$ (Alternatively, one can
say that $\psi^{-1}Y = \tilde{E}.$) Let the relative canonical
divisor of $\psi$ be $\sum_{i \in J}(\nu_i-1)E_i,$ i.e. $\nu_i-1$
is the multiplicity of $E_i$ in the divisor of $\psi^*(dx_1 \wedge
\ldots \wedge dx_n).$ Finally put $E_I^\circ := (\cap_{i \in I}
E_i) \< (\cup_{l \not \in I}E_l)$ for $I \subset J;$ these
$E_I^\circ$ form a natural locally closed stratification of
$\tilde{X}.$ (Note that $E_{\emptyset}^\circ = \tilde{X} \<
\cup_{l \in J}E_l.$) The (local) \emph{topological zeta function
of $\I$ at $0$} is
$$Z_{top, \I}(s) := \sum_{I\subset J} \chi(E_I^\circ \cap
\psi^{-1}\{0\})\prod_{i \in I} \frac1{\nu_i + s N_i} \in
\mathbb{Q}(s).$$ There is a global version replacing $E_I^\circ
\cap \psi^{-1}\{0\}$ by $E_I^\circ.$

When $r=1,$ Denef and Loeser showed in \cite{DL-zeta fie onafh van
res} that the expression above does not depend on the chosen
principalization (which for $r=1$ is just an embedded resolution)
by writing it as a limit of $p$-adic Igusa zeta functions.
Alternatively, they obtained it later in \cite{DL-motivische
zetafuncties} as a specialization of the motivic zeta function.
This can be generalized to arbitrary $r,$ see e.g.
\cite[(2.4)]{veys-zuniga}. Still another possibility for arbitrary
$r$ is to use the Weak Factorization Theorem of W\l odarczyk et
al. \cite{AKMW-factorization} to compare two principalizations. At
any rate, observe that a complete list of possible poles of
$Z_{top, \I}(s)$ is given by the $-\frac{\nu_i}{N_i}, i \in J.$

\vspace{\baselineskip}

For $r=1,$ say $\I = (f),$ there are the following intriguing
conjectures \cite{DL-zeta fie onafh van res}. (They were
originally formulated for the $p$-adic Igusa zeta function, which
is a certain $p$-adic integral, and partially motivated by
analogous statements that are true for a similar complex integral.
See \cite{Igusa}, in particular section 5.4.)

\begin{description}
\item{\underline{Conjecture 1.}} If $s_0$ is a pole of $Z_{top,
f}(s),$ then $s_0$ is a root of the (local) Bernstein-Sato
polynomial $b_{f,0}(s)$ of $f.$

\item{\underline{Conjecture 2.}} \emph{(Monodromy Conjecture.)} If
$s_0$ is a pole of $Z_{top, f}(s),$ then $e^{2\pi i s_0}$ is an
eigenvalue of the local monodromy action on some cohomology group
of the Milnor fibre of $f$ at some point of $\{f=0\}$ close to
$0.$
\end{description}
Note that Conjecture 1 implies Conjecture 2 since for any root
$s_0$ of $b_{f,0}(s)$ we have that  $e^{2\pi i s_0}$ is such a
monodromy eigenvalue \cite{Malgrange - Bernstein}. For $n=2$
Conjecture 1 was proved by Loeser \cite{Loeser num data}. He also
verified it for non-degenerate polynomials satisfying extra
assumptions \cite{Loeser 2}. (There is a more elementary proof of
Conjecture 2 for $n=2$ by Rodrigues \cite{Rodrigues monodromie}.)
Concerning Conjecture 2, there are various partial results, mainly
for $n=3,$ by Artal, Cassou-Nogu\`es, Luengo and Melle
\cite{ACLM1}, \cite{ACLM2}, and Lemahieu, Rodrigues and the second
author \cite{veys 1993 poles and monodromy}, \cite{Rodrigues-Veys
monodromie in dim3}, \cite{veys 2006 vanishing},
\cite{lemahieu-veys monodromie op oppervlakken}, \cite{LV2}.

\vspace{\baselineskip}

Budur, Musta\c{t}\v{a} and Saito introduced a Bernstein-Sato
polynomial associated to an arbitrary ideal $\I = (f_1, \ldots ,
f_r) \subset \C[x_1, \ldots , x_n].$ (Their polynomial coincides
with a polynomial that appears in \cite{Gyoja}.)


Still for arbitrary $r,$ the notion of local Milnor fibre is in
general not well-defined. There is however the following
construction of Verdier. To any constructible complex of sheaves
$\mathcal{F}^\bullet$ on $\mathbb{A}^n,$ he associates a similar
complex on $C_Y\mathbb{A}^n,$ the normal cone of $Y =
\mathbf{V}(\I)$ in $\mathbb{A}^n.$ This complex is called the
\emph{specialization} of $\mathcal{F}^\bullet.$ It is moreover
equipped with a `canonical monodromy operator'. In particular when
$\I = (f)$ it turns out that this specialization of $\C^\bullet$
is in some sense equivalent to the usual complex of nearby cycles
on $Y = \mathbf{V}(f)$, and that the two monodromy notions
essentially correspond.

A nice feature of the Bernstein-Sato polynomial of \cite{Budur -
Mustata - Saito} is that for any of its roots $s_0,$ we have that
$e^{2 \pi i s_0}$ is a `Verdier monodromy eigenvalue', see
\cite[Corollary 2.8]{Budur - Mustata - Saito}. This thus
generalizes the implication for $r=1$ mentioned above.

\vspace{\baselineskip}

It is natural to ask also for arbitrary $\I$ if poles of the
topological zeta function of $\I$ are always roots of its
Bernstein-Sato polynomial, remembering of course that this
question turned out to be very difficult already for $r=1.$ In the
special case of a monomial ideal $\I$ this was verified in
\cite{mustata & co: bernstein voor monomiale idealen} by Howald,
Musta\c t\v a and Yuen (for the $p$-adic Igusa zeta function).
Maybe more accessible, do the poles of the topological zeta
function of an arbitrary ideal $\I$ induce monodromy eigenvalues
in the sense of Verdier? The main result of this paper is to
provide an affirmative answer to the last question for arbitrary
ideals in two variables.

\vspace{\baselineskip}

The plan of the paper is as follows. We work over the base field
of complex numbers. In \S 2 we explain the construction of the
specialization functor of Verdier. In arbitrary dimension we show
in \S 3 a formula for the `Verdier monodromy eigenvalues' of an
ideal $\I$ in terms of a principalization of the ideal, in the
same spirit as A'Campo's formula \cite{A'Campo} for the
eigenvalues of one function $f$ in terms of an embedded resolution
of $Y = \mathbf{V}(f).$ Note that this is a priori not obvious;
for $r=1$ the complex of nearby cycles lives on $Y$ but for $r >
1$ the specialized complex of Verdier lives not on $Y$ but on
$C_Y\mathbb{A}^n.$ In \S 4 we prove the `generalized Monodromy
Conjecture' for ideals in two variables, and finally in \S 5 we
provide some examples.

\vspace{\baselineskip}

\noindent \textsc{Acknowledgements:} The authors want to thank J.
Sch\"urmann for the interesting conversations and explanations.

\section{The specialization functor of Verdier} \label{sectie Verdier}

Let $\I$ be a coherent ideal sheaf on a  variety $X.$ Consider the
associated subscheme $Y = \mathbf{V}(\I)$ of $X .$ We construct
the blow-up $\pi : B = Bl_Y\,X \to X$ of $X$ in $Y$ and denote by
$E$ the inverse image $\pi^{-1}(Y).$

For every $e \in E$ we want to study the zeta function of
monodromy of the ideal $\I.$ To define what monodromy is in this
context, we will need notions as the normal cone of $Y$ in $X$ and
the specialization functor of Verdier.

The normal cone of $Y$ in $X$ is a cone over $Y$ defined as
$$C_YX = \mbox{Spec}_{\mathcal{O}_Y}(\oplus_{n \geq 0} \I^n / \I^{n+1})\, ,$$
see e.g. \cite[B.6]{Fulton-intersection theory}. Interesting to
notice is that the projectivization $P(C_YX)$ of the normal cone
is exactly the exceptional variety $E$ of the blowing-up of $X$ in
$Y$ (with its non-reduced scheme structure). We can identify the
`locus of vertices' of this cone with $Y$. So we have an embedding
$j : Y \to C_YX$ and a projection $p: C_YX \to Y.$ These maps
satisfy $p \circ j = Id_Y.$ Moreover, we have an action of $\C^*$
on the normal cone, coming from the graduation on $\oplus_{n \geq
0} \I^n / \I^{n+1}.$ The locus of vertices is the scheme of the
fixed points of the action and the morphisms $j$ and $p$ commute
with it. Starting with a point $e \in (C_YX) \< Y,$ we define the
\emph{ruler} through $e$ as the orbit of $e$ under the action of
$\C^*$ (with its reduced scheme structure), so every ruler can be
identified with $\C^*.$

The construction of the normal cone is functorial in the following
sense. Suppose we have a map $f: X' \to X$ of schemes and two
subschemes $Y \subset X$ and $Y' \subset X'$ such that
$$
\begin{array}{ccc} Y' & \hookrightarrow & X'
 \\
\ \  \ \ \downarrow \, \scriptstyle{f|_{Y'}} & & \ \  \downarrow \, \scriptstyle{f}  \\
Y & \hookrightarrow & X
\end{array}$$
is a cartesian diagram, then we can associate a map $C(f):
C_{Y'}X' \to C_YX$ to $f.$ In particular, the map $\pi : B \to X$
has an associated map $C(\pi) : C_EB \to C_YX.$ This last one has
the following interesting properties which we will use further on.
If we restrict the map $C(\pi)$ to the `punctured' cones, $C'(\pi)
: (C_EB) \< E \to (C_YX) \< Y,$ this is an isomorphism. Likewise,
the map of projectivized cones $C''(\pi) : P(C_EB) \to P(C_YX)$ is
well-defined and is an isomorphism. Moreover, the following
diagram is commutative.

\begin{equation} \label{comm diagram C(pi) Schurmann}
\begin{array}{ccccc} E & \stackrel{\sim}{\longleftarrow} & P(C_EB) & \longleftarrow
& (C_EB)\<E \\
& & & & \\

 \ \ \Big\downarrow \, \scriptstyle{\pi} &  \circ & \ \ \wr \!
\Big\downarrow \, \scriptstyle{C''(\pi)}
  & \circ &  \ \wr \!  \Big\downarrow \, \scriptstyle{C'(\pi)}\\
  & & & & \\
Y & \longleftarrow & P(C_YX) & \longleftarrow & (C_YX)\< Y
\end{array}
\end{equation}

The (canonical) deformation to the normal cone of $Y$ in $X$ is a
scheme $\hat{X}$ with a morphism $\hat{\pi} : \hat{X} \to
\mathbb{A}^1$ such that $\hat{\pi}^{-1}(\mathbb{A}^1 \setminus
\{0\}) \cong X \times (\mathbb{A}^1 \setminus \{0\})$ and the
fibre over $0 \in \mathbb{A}^1$ is the normal cone $C_Y X.$ This
construction can be found in \cite[Chapter 5]{Fulton-intersection
theory}. So we obtain the following commutative diagram.
$$
\begin{array}{ccccc} C_YX & \longrightarrow & \hat{X}
& \stackrel{j}{\hookleftarrow} &
X \times (\mathbb{A}^1 \setminus \{0\}) \\
\Big\downarrow & \circ  & \ \ \Big\downarrow  \scriptstyle{\hat{\pi}} & \ \ \ \circ  & \Big\downarrow  \\
0 & \longrightarrow & \mathbb{A}^1 &
\stackrel{j_0}{\hookleftarrow} & \mathbb{A}^1 \setminus \{0\}
\end{array}$$

We introduce now the specialization functor of Verdier that
associates (by means of the previous diagram) to a constructible
complex of sheaves $\mathcal{F}^{\bullet}$ on $X$ a constructible
complex of sheaves on $C_YX$.

\bigskip
\noindent \emph{Note.} By a constructible complex of sheaves (of
$\C$-vector spaces) $\mathcal{G}^\bullet$ on a scheme $Z$ we mean
an object in the full subcategory $D^b_c(Z)$ of the derived
category $D(Z)$; its cohomology sheaves $H^n(\mathcal{G}^\bullet)$
are thus constructible and moreover zero for $|n|>>0$. (A sheaf on
$Z$ is called constructible if there is a Zariski-locally closed
stratification of $Z_{red}$ such that the restriction of the sheaf
to each stratum is locally constant for the complex topology.)

\bigskip
In \cite{Verdier artikel} the specialization of
$\mathcal{F}^{\bullet}$ is defined as
$$Sp_{Y \backslash X}(\mathcal{F}^{\bullet}) :=
\psi_{\hat{\pi}}(j_! \, pr_1^* \, \mathcal{F}^{\bullet}),$$ where
$\hat{\pi}$ and $j$ are as before, $pr_1$ is the projection of $X
\times (\mathbb{A}^1 \setminus \{0\})$ on the first factor and
$\psi_{\hat{\pi}}$ is the nearby cycle functor of Deligne
\cite{Deligne-monodromie met vanishing cycles} (see also
\cite[section 4.2]{Dimca}). In \cite[p.356-357]{Verdier artikel},
Verdier defines a canonical transformation of monodromy on the
complex of sheaves $Sp_{Y \backslash X} (\mathcal{F}^{\bullet}).$

This specialization functor has a number of important properties,
from which we will give the ones that we need in this article.
These properties are stated in \cite[sections 8 and 10]{Verdier
artikel}.

\begin{description}
\item{(SP1) \emph{Monodromy.}} For every constructible complex of
sheaves $\mathcal{F}^\bullet,$ the complex $Sp_{Y \<
X}(\mathcal{F}^\bullet)$ is monodromic. This means that $Sp_{Y \<
X}(\mathcal{F}^\bullet)$ is locally constant (with respect to the
usual complex topology) on every ruler of $C_YX.$

\item{(SP2) \emph{Proper direct image.}} Suppose we have a
cartesian diagram
$$
\begin{array}{ccc} Y' & \hookrightarrow & X'
 \\
\downarrow & & \ \  \downarrow \, \scriptstyle{f}  \\
Y & \hookrightarrow & X
\end{array}$$
where $f$ is proper. Let $C(f): C_{Y'}X' \to C_YX$ be the morphism
associated to $f.$ Then $C(f)$ is proper and for each
constructible complex $\mathcal{F}^\bullet$ on $X',$ the natural
morphism
$$ Sp_{Y \< X}(Rf_* \mathcal{F}^\bullet) \to RC(f)_*(Sp_{Y' \<
X'}(\mathcal{F}^\bullet))$$ is an isomorphism.

\item{(SP6) \emph{Normalization.}} Suppose that $Y$ is a principal
divisor with equation $f=0.$ Then the morphism $p \times C(f):
C_YX \to Y \times \mathbb{A}^1$ is an isomorphism, where $p: C_YX
\to Y $ is the projection. Notice that $C_{\{0\}}\mathbb{A}^1 =
\mathbb{A}^1,$ so $C(f)$ is indeed a map $C_YX \to \mathbb{A}^1.$
With this isomorphism, we can define a section for $z \in \C \<
\{0\}:$
\begin{eqnarray*}
s_z \  :\ \  Y & \to & C_YX \\
y & \mapsto & (p \times C(f))^{-1}(y,z).
\end{eqnarray*}
Then we have for each constructible complex of sheaves
$\mathcal{F}^\bullet$ on $X$ an isomorphism
$$s_z^* (Sp_{Y\< X}(\mathcal{F}^\bullet)) \to
\psi_f(\mathcal{F}^\bullet).$$ This isomorphism is compatible with
the two monodromy operations defined on both sheaves on $Y,$ in
the sense that they are each other's opposite. (We refer to
\cite{Deligne-monodromie met vanishing cycles} or \cite[section
4.2]{Dimca} for a description of the \lq classical\rq\ monodromy
on $\psi_f(\mathcal{F}^\bullet)$.)

\item{(SP7) \emph{Perversity.}} The specialization functor
transforms perverse sheaves into perverse sheaves.

\end{description}

\section{The zeta function of monodromy}

As we said before, Verdier defined a canonical transformation of
monodromy on the complex of sheaves $Sp_{Y\< X}(\C_X^\bullet),$
which we will denote by $M.$ For each $m \in \Z$ and $y \in (C_YX)
\<Y,$ we have an automorphism $M_y^m$ on the stalk
$H^m(Sp_{Y\<X}(\C_X^\bullet)_y).$ \emph{Eigenvalues of monodromy}
are eigenvalues of these vector space automorphisms.

We now define the zeta function of monodromy $Z_{\I,e} (t)$ of
$\I$ for a point $e \in E$. Choose an arbitrary point $e'$ on
$(C_YX) \backslash Y$ that is mapped to $e$ by the
projectivization. Then $Z_{\I,e} (t)$ is the (finite) product
$$Z(Sp_{Y\< X}(\C_X^{\bullet}))(e'):=
\prod_{m \in \Z} \mbox{det} (Id - t \,
M^m_{e'})^{(-1)^m}\, .
$$
Note that this is the usual notion of a zeta function for $M$ and
$e'$, and that, by (SP1), this definition is independent of the
choice of $e'$.

In the next section, we use the zeta function of monodromy to
prove the monodromy conjecture in dimension two. We will prove
that a pole $s_0$ of the topological zeta function induces a zero
or a pole $e^{2\pi i s_0}$ of the monodromy zeta function for some
point $e \in E.$ This implies that this number $e^{2\pi i s_0}$ is
an eigenvalue of monodromy.

\begin{remark} Note that it is also true that each eigenvalue of
monodromy is a zero or a pole of the monodromy zeta function for
some point $e \in E.$ Since we know from (SP7) that $Sp_{Y
\backslash X} (\C_X^{\bullet})$ is a (shifted) perverse sheaf, we
can copy the proof of \cite[Lemma 4.6]{Denef - monodromie} to
obtain this result.
\end{remark}

Now we prove  a generalization of the formula of A'Campo
\cite[Theorem 3]{A'Campo} to the case of ideals.

\begin{theorem}
Let $\I$ be a sheaf of ideals on a variety $X.$ Let $Y
=\mathbf{V}(\I)$ be the associated subscheme of $X$ and suppose
that $\mbox{Sing}\, (X) \subset Y.$ Let $\pi : B = Bl_Y\,X \to X$
be the blow-up of $X$ in $Y$ and   $\psi : \tilde{X} \to X$ a
principalization of $\I.$ Define $\fie : \tilde{X} \to B$ as the
unique morphism such that $\psi = \pi \circ \varphi.$ We denote by
$E$ the inverse image $\pi^{-1}(Y)$ and $\tilde{E} =
\psi^{-1}(Y).$ We use $E_i,$ $i \in J,$ for the irreducible
components of $\tilde{E}$ and $N_i$ for the according
multiplicities in $\tilde{E}.$ Put $E_i^\circ = E_i \setminus
\cup_{j \in J ,\, j\neq i} E_j.$ For a point $e \in E,$ the zeta
function of monodromy is equal to
$$Z_{\I,e} (t) =
\prod_{j \in J} (1-t^{N_j})^{\chi(E_j^\circ \cap
\varphi^{-1}(e))}.$$
\end{theorem}

\begin{remark} (1) Note that the principalization $\psi$ factorizes
through the blow-up $\pi$ since the inverse image $\psi^{-1}(Y)$
is a Cartier subscheme in $\tilde{X}.$

(2) When $\I$ is principal we can consider $\pi$ as the identity
and thus $\psi=\varphi$. Then $Z_{\I,e} (t)$ is the usual
hypersurface zeta function of monodromy for $e\in E=Y$ and we
recover A'Campo's formula.
\end{remark}

\noindent \emph{Proof.} Fix $e \in E$ and choose an element $e'
\in (C_YX) \backslash Y$ that is mapped to $e$ by the
projectivization map.

The restriction of $\pi$ is an isomorphism $B \setminus E \to X
\setminus Y,$ which implies that $R\pi_* \C_B^\bullet |_{X
\setminus Y} = \C_X^\bullet |_{X \setminus Y}.$ Since $ Sp_{Y
\backslash X}(\mathcal{F}^\bullet) |_{(C_YX)\setminus Y}$ only
depends on $\mathcal{F}|_{X \setminus Y}$ (see
\cite[p.354]{Verdier artikel}), we can deduce that
$$Sp_{Y \backslash X}(\C_X^\bullet) |_{(C_YX)\setminus Y} = Sp_{Y \backslash
X}(R \pi_* \C_B^\bullet)|_{(C_YX)\setminus Y}$$ and $$Z(Sp_{Y
\backslash X}(\C_X^\bullet))(e')  = Z(Sp_{Y \backslash X}(R \pi_*
\C_B^\bullet))(e').$$

Because we have the cartesian diagram
$$
\begin{array}{ccc} E & \hookrightarrow & B
 \\
\downarrow & & \ \ \ \downarrow \pi  \\
Y & \hookrightarrow & X
\end{array}$$
we can use (SP2) to write that
$$Sp_{Y \backslash X}(R \pi_* \C_B^\bullet) \stackrel{\sim}{\longrightarrow}
RC(\pi)_*(Sp_{E \backslash B}(\C_B^\bullet))$$ is an isomorphism,
or that
$$ Z(Sp_{Y \backslash X}(R \pi_* \C_B^\bullet))(e') =
Z(Sp_{E \backslash B}(\C_B^\bullet))(C(\pi)^{-1}(e')).$$ (Note
that the restriction of $C(\pi)$ to $(C_EB) \< E$ is an
isomorphism, as said in \S \ref{sectie Verdier}.)

There exists an open $V_e$ around $e$ such that $E \cap V_e$ is
the zero locus of one nonzerodivisor $f.$ So we can use (SP6) to
see that locally, we have an isomorphism $p \times C(f): C_EB \to
E \times \mathbb{A}^1.$ Take $z \in \C \< \{0\}$ such that $(p
\times C(f))(C(\pi)^{-1}(e')) = (e,z).$ This is possible, using
the commutative diagram (\ref{comm diagram C(pi) Schurmann}). We
can conclude that we have an isomorphism
$$s_z^*(Sp_{E \< B}(\C_B^\bullet)) \longrightarrow \psi_f(\C_B^\bullet).$$
This isomorphism is `compatible' with the monodromy actions, so
$$Z(Sp_{E\< B}(\C_B^\bullet))(s_z(e)) =
Z(\psi_f(\C_B^\bullet))(e).$$ The right hand side is the monodromy
zeta function of the map $f: V_e \subset B \to \C$ in the point
$e.$ Because we already have an embedded resolution of $E$ in $B,$
namely the map $\fie: \tilde{X} \to B,$ we can use the formula of
A'Campo \cite[Theorem 3]{A'Campo}. (This formula was originally
proven for functions on a smooth variety. But one can check that
for instance the proof that Dimca gives in \cite[Corollary
6.1.15]{Dimca} is still valid in this more general context where
the hypersurface contains the singular locus of the ambient
variety.) We obtain
$$ Z(\psi_f(\C_B^\bullet))(e) = \prod_{j \in J} (1-t^{N_j})^{\chi(E_j^\circ
\cap \varphi^{-1}(e))}.$$ For this last equality, we need that the
multiplicities of the irreducible components $E_i$ of $\tilde{E}$
in the divisor $\fie^{-1}(E)$ are precisely the $N_i$. Indeed:
$$\tilde{E}=\sum_{i \in J} N_iE_i = (\pi \circ \fie)^{-1}(Y) = \fie^{-1}(\pi^{-1}Y) =
\fie^{-1}(E).$$ Putting all these equalities of zeta functions
together, proves our theorem. \hfill$\square$

\section{The Monodromy Conjecture}

From now on, we will work in dimension 2, but we will use the same
notation as before. So let $\I \subset \C[x,y]$ be an ideal
satisfying $\{0\} \subset \mbox{Supp} \, \I.$ Put $Y = V(\I)$ the
subscheme of $X = \C^2$ defined by $\I.$ Let $\pi: B=Bl_YX \to X$
be the blowing-up of $X$ in $Y$ and $\psi: \tilde{X} \to X$ be the
minimal principalization of $\I.$ The map $\fie: \tilde{X} \to B$
is defined such that $\psi = \pi \circ \fie.$ We denote
$E=\pi^{-1}(Y)$ as divisor on $B$, and $\tilde{E} = \psi^{-1}(Y) =
\sum_{i \in J}N_iE_i$ as divisor on $\tilde{X}$. Here the
$E_i(N_i, \nu_i)$ are the irreducible components of $\tilde{E}$,
together with their numerical data.

\begin{theorem} \emph{(Generalized Monodromy Conjecture.)} \label{monodromieconjectuur idealen}
    If $-\frac{\nu}N$ is a pole of the local topological zeta
    function of an ideal $\I \subset \C[x,y],$ then there exists a point $y \in E$
    such that $e^{-2\pi i \frac{\nu}{N}}$ is an eigenvalue of
    monodromy in $y.$
\end{theorem}

We  first mention some results that will be useful in the proof.

\vspace{\baselineskip}

For an ideal $\I = (f_1, \ldots , f_r) \subset \C[x,y]$ we can
look at the linear system $\{\lambda_1f_1 + \ldots + \lambda_rf_r
= 0 \, | \, \lambda_i \in \mathbb{C} \mbox{ for } i =1, \ldots,r
\}.$ A \emph{total generic curve} of $\I$ is a general element of
this linear system. Now we determine whether there are common
components among the $f_i$ and put them together. So we can write
$\I = (h)(f_1', \ldots , f_r'),$ where $(f_1', \ldots , f_r')$ is
a finitely supported ideal. A \emph{generic curve} of the ideal
$\I$ is a general element of the linear system $\{\lambda_1f_1' +
\ldots + \lambda_rf_r' = 0 \, | \, \lambda_i \in \mathbb{C} \mbox{
for } i =1, \ldots,r \}.$ Notice that the definition of the
(total) generic curve depends on the choice of generators we use
to represent the ideal.

\begin{lemma} \label{lemma gecontracteerd -> geen snijding}
If an irreducible component $E_i$ of $\tilde{E}$ is contracted by
$\varphi,$ then there is no intersection between $E_i$ and the
strict transform of a generic curve of the ideal $\I$ in
$\tilde{X}.$
\end{lemma}
\noindent \emph{Proof.} If there is an intersection, the strict
transform of a generic curve of $\I$ in $B$ contains the point
$\fie(E_i).$ But if we denote $\widetilde{f_i}$ for the strict
transform of $f_i'$ in $B,$ the linear system
$\{\lambda_1\widetilde{f_1} + \ldots + \lambda_r\widetilde{f_r}=0
\, | \, \lambda_i \in \mathbb{C} \mbox{ for } i =1, \ldots,r \}$
should be base point free (see e.g. \cite[Example
7.17.3]{Hartshorne}). \hfill$\square$

\begin{lemma}\label{lemma alpha-relatie idealen}
Let $(f_1, \ldots , f_l)$ be an ideal in $\mathbb{C}[x,y]$ and
$E_0(N, \nu)$  an exceptional curve of the principalization
(together with the numerical data). Suppose $E_0$ intersects $n$
 times the strict transform of a  generic curve and $m$ times other components
 $E_1
(N_1, \nu_1), \ldots , E_m (N_m, \nu_m)$  of the principalization.
Put $\alpha_i = \nu_i - \frac{\nu}{N}N_i$ for $i=1, \ldots, m.$
Then we have
\begin{enumerate}
\item $\sum^{m}_{i=1} \alpha_i = m-2+\frac{\nu n}{N},$ where the
left hand side is zero when $m=0$, and \item $-1 \leq \alpha_i <
1$ for $i \in \{1, \ldots , m\}$. Moreover, $\alpha_i=-1$ only
occurs when $m=1$.
\end{enumerate}
\end{lemma}

The first equality can be found in \cite[$\S 3$]{polen ideaal
dim2}. It is a reformulation of the relation between the numerical
data proved by Loeser in \cite{Loeser num data}. The second
statement is \cite[Proposition 3.1]{polen ideaal dim2}.

 \vspace{\baselineskip}

In \cite{polen ideaal dim2}, we gave a complete list of five
conditions in which a candidate pole is indeed a pole of the local
topological zeta function of an ideal in dimension two. Using
Lemma \ref{lemma alpha-relatie idealen}, it is not difficult to
see that the conditions given in numbers 2, 3 and 4 of
\cite[Theorem 4.2]{polen ideaal dim2} are equivalent to the
condition of being intersected by the strict transform of a
generic curve. So we can reformulate this theorem as follows.

\begin{theorem} \label{stelling polen zeta functie ideaal 2} Let $\mathcal{I} \subset \mathbb{C}[x,y]$ be an ideal
satisfying $0 \in \mbox{Supp}\,(\mathcal{I})$ and $\psi: \tilde{X}
\to \mathbb{C}^2$ the minimal principalization of $\mathcal{I}$ in
a neighbourhood of $0.$ Let $E_\bullet(N_\bullet,\nu_\bullet)$ be
the components of the support of the total transform
$\psi^*\mathcal{I}$ with their associated numerical data.\\
The rational number $s_0$ is a pole of the local topological zeta
function of $\mathcal{I}$ if and only if  at least one of the
following conditions is satisfied:
\begin{enumerate}
\item $s_0 = -\frac1{N}$ for a component $E_0(N, \nu)$ of the
support of the weak transform of $\mathcal{I}$; \item $s_0 =
-\frac{\nu}{N}$ for $E_0(N, \nu)$ an exceptional curve that has
non-empty intersection with the strict transform of a generic
curve of the ideal $\I;$ \item $s_0 = -\frac{\nu}{N}$ for $E_0(N,
\nu)$ an exceptional curve that intersects at least three times
other components.
\end{enumerate}
\end{theorem}

Now we are ready to prove the monodromy conjecture for ideals in
dimension 2. We were inspired by \cite{Rodrigues monodromie},
where Rodrigues gives an elementary proof of the monodromy
conjecture for curves on normal surfaces. \vspace{\baselineskip}

\noindent \emph{Proof of Theorem \ref{monodromieconjectuur
idealen}.} Choose a pole $s_0$ of the local topological zeta
function of the ideal $\I$ and take $a,d \in \Z$ satisfying $s_0 =
-\frac{a}{d}$ and $\mbox{gcd}(a,d) = 1.$

Suppose  that there exists an irreducible component $E_m (N_m ,
\nu_m)$ with $d\mid N_m$ that is not contracted to a point by
$\fie: \tilde{X} \to B.$ Choose a point $y \in \varphi(E_m^\circ)$
such that $\varphi^{-1}(y)$ is a finite set of points. Then there
exists for every $j \in J$ a nonnegative integer $k_j$ (with $k_m
\neq 0$) such that
$$Z_{\I,y} (t) = \prod_{j \in J} (1-t^{N_j})^{k_j},$$
from which we see that $e^{2\pi i s_0}$ is a zero of the zeta
function of monodromy in $y,$ so it is an eigenvalue of monodromy.

On the other hand, suppose that every $E_i(N_i, \nu_i)$ that
satisfies $d\mid N_i$  is  contracted to a point by $\fie.$ This
implies (by using Lemma \ref{lemma gecontracteerd -> geen
snijding} and Theorem \ref{stelling polen zeta functie ideaal 2})
that there exists an index $m \in J,$ with $s_0 =
-\frac{\nu_m}{N_m},$ such that $E_m$ is exceptional and $E_m$
intersects at least three times other components of $\tilde{E}.$
Fix such a $m$ and take $y = \fie(E_m).$  Define $T \subset J$
such that $\fie^{-1}\{y\} = \cup_{i \in T} E_i.$ To prove that
$e^{2 \pi i s_0 }$ is a zero or a pole of the zeta function of
monodromy in $y,$ it is enough to show that
$$\sum_{i \in T ,\, d \mid  N_i} \chi(E_i^\circ) < 0.$$

First notice that the strict transform of a generic curve
intersects at least one exceptional component, so not every
exceptional curve can be contracted and it is impossible that $d$
divides $N_i$ for every $i\in J.$ We define
$$\tilde{M} := \bigcup_{i \in T , \, d \mid N_i} E_i.$$

Every connected component of $\tilde{M}$ contains at least one
irreducible component that intersects a component $E_k(N_k,
\nu_k)$ outside $\tilde{M}$. If $d \mid N_k$, then $E_k$ is
contracted by $\fie$ and since $E_k$ has a non-empty intersection
with $\tilde{M},$ this also implies that $k \in T$. So $d \nmid
N_k$.
Moreover, the irreducible component of $\tilde{M}$ that intersects
$E_k$ needs to intersect a second component $E_{k'}(N_{k'},
\nu_{k'})$ with $d \nmid N_{k'}.$ We can deduce this from the
formula $\kappa N_l = \sum_{i=1}^r N_i,$ where $E_l$ is an
exceptional curve intersecting $r$ times other components $E_i,$
for $i = 1, \ldots , r,$ and where $-\kappa$ denotes the
self-intersection number of $E_l.$

Every connected component $\bigcup_{i \in I_0} E_i$ of $\tilde{M}$
satisfies
$$\sum_{i \in I_0} \chi(E_i^\circ) \leq 0.$$
Because such a component has at least two external intersections,
this is a direct consequence of the following fact: if
$\bigcup_{i=1}^rE_i$ is a tree consisting of rational curves, then
$\sum_{i=1}^{r} \chi(E_i \< (\cup_{j = 1 , \ldots , r , \, j\neq
i}E_j)) = 2$. See e.g. \cite[Lemma 2.2]{Rodrigues monodromie}.

Now fix the connected component $M = \bigcup_{i \in I_M} E_i $ of
$\tilde{M}$ that contains $E_m.$ For this one, we prove the strict
inequality $$\sum_{i \in I_M} \chi(E_i^\circ) < 0.$$ If this sum
 would equal 0, then $M$ is a tree of rational curves that
intersects precisely two times with components $E_j$ outside
$\tilde{M}.$ From a previous argument, we know that these
intersections are on one curve $E_s$ with $s \in I_M.$ There are
two possibilities.
\begin{description}
\item{$E_s = E_m.$} The component $E_m$ will intersect at least
 one other curve $E_j(N_j, \nu_j).$ This curve belongs to $M,$
so $d \mid N_j$ and the number $\nu_j - \frac{a}{d} N_j$ is an
integer. The only possibility is 0 (see Lemma \ref{lemma
alpha-relatie idealen}), so this implies that $-\frac{\nu_j}{N_j}
= -\frac{a}{d}.$ We know that $E_j$ is contracted to a point, so
it is not intersected by the strict transform of a generic curve.
From the relations between the numerical data (see also Lemma
\ref{lemma alpha-relatie idealen}), we can deduce that $E_j$
intersects at least one other component of $M$ and we can repeat
this process infinitely many times. This leads to a contradiction.

\item{$E_s \neq E_m.$} Now we can start for each of the (at least
three) components that intersect $E_m$ the same procedure of
constructing a series of curves in $M.$ Only one of these can
eventually stop, because there is only one of these series that
has intersections outside of $M.$
\end{description}
This ends the proof. \hfill$\square$

\begin{remark}
Theorem \ref{monodromieconjectuur idealen} is also true for the
motivic and the Hodge zeta function of $\I,$ and for $p$-adic
Igusa zeta functions associated to several polynomials in two
variables, because the necessary condition of \cite[Theorem
4.2]{polen ideaal dim2} to be a pole is still valid for these zeta
functions. We refer to \cite[$\S 6$]{polen ideaal dim2} for
definitions and more explanation.

\end{remark}

\section{Some examples}

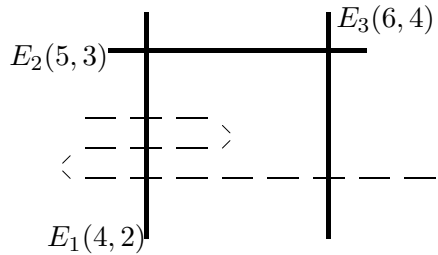
\begin{figure}[h]  \caption{\textsl{Intersection diagram of the principalization
of $(x^3y, x^6 + y^4).$}}\label{tekening voorbeeld 1}
\unitlength=1mm
\begin{center}
\begin{picture}(60,28)(10,10)
\linethickness{0.5mm}

\put(32,7){\line(0,1){30}} \put(27,32){\line(1,0){34}}
\put(56,7){\line(0,1){30}}

\put(19,6){$E_1(4,2)$} \put(14,30){$E_2(5,3)$}
\put(57,35){$E_3(6,4)$}

\linethickness{0.15mm}

\put(24,23){\line(1,0){4}} \put(30,23){\line(1,0){4}}
\put(36,23){\line(1,0){4}}

\put(42,22){\line(1,-1){1}} \put(42,19){\line(1,1){1}}
\put(24,19){\line(1,0){4}} \put(30,19){\line(1,0){4}}
\put(36,19){\line(1,0){4}} \put(22,18){\line(-1,-1){1}}
\put(22,15){\line(-1,1){1}} \put(24,15){\line(1,0){4}}
\put(30,15){\line(1,0){4}} \put(36,15){\line(1,0){4}}
\put(42,15){\line(1,0){4}}

\put(48,15){\line(1,0){4}} \put(54,15){\line(1,0){4}}
\put(60,15){\line(1,0){4}} \put(66,15){\line(1,0){4}}
\end{picture}
\end{center}
\end{figure}

(1)  $\I = (x^3y, x^6 + y^4) \subset \C[x,y]$. The intersection
diagram of the principalization $\tilde{X}$ is as in Figure
\ref{tekening voorbeeld 1}, where the dashed line denotes the
strict transform of a generic curve.

\noindent We see that the poles of the local topological zeta
function are $-\frac12$ and $-\frac23.$

The blowing-up $B = Bl_Y\C^2$ is given by
$$\mbox{Proj} \frac{\C[x,y][A,B]}{(Ax^3y-B(x^6+y^4))}, $$
see e.g. \cite[Section IV.2.1]{geometry of schemes}. The
exceptional variety consists of only one irreducible curve, say
$E.$ The map $\fie : \tilde{X} \to B$ contracts $E_2$ to a point
$a$ on $E.$ It maps $E_1$ and $E_3$ surjectively on $E.$ The
restriction to $E_3$ is one-to-one, but the restriction to $E_1$
is three-to-one.

For a point $e$ on $E$ different from $a,$ we compute that the
monodromy zeta function is
$$Z_{\I,e}(t)=(1-t^4)^3(1-t^6).$$
We can immediately check that $e^{-\pi i}$ and $e^{-\frac{4 \pi
i}{3}}$ are zeroes of this function. Since $-\frac12$ and
$-\frac23$ were the only poles of the local topological zeta
function, this illustrates the monodromy conjecture.

\vspace{\baselineskip}

\begin{figure}[h]  \caption{\textsl{Intersection diagram of the principalization
of $(x^4, xy^2, y^3) .$}}\label{principalisatie voorbeeld 2}
\unitlength=1mm
\begin{center}
\begin{picture}(60,28)(10,10)
\linethickness{0.5mm}

\put(32,12){\line(0,1){25}} \put(27,32){\line(1,0){34}}
\put(56,12){\line(0,1){25}}

\put(19,11){$E_1(3,2)$} \put(14,30){$E_3(8,5)$}
\put(57,35){$E_2(4,3)$}

\linethickness{0.15mm}

\put(25,22){\line(1,0){3}} \put(30,22){\line(1,0){4}}
\put(36,22){\line(1,0){3}}

\put(42,24){\line(-1,-1){2}}

\put(43,25){\line(0,1){3}} \put(43,30){\line(0,1){4}}
\put(43,36){\line(0,1){3}}

\end{picture}
\end{center}
\end{figure}
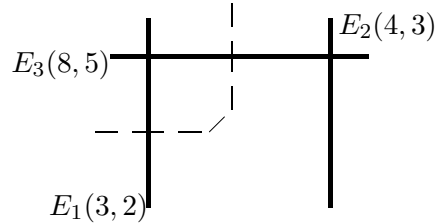

(2) $\I = (x^4, xy^2, y^3) \subset \C[x,y]$. The intersection
diagram of the principalization (with numerical data) is given in
Figure \ref{principalisatie voorbeeld 2}. We can easily deduce
that the poles of $Z_{top, \I}(s)$ are $-\frac23$ and $-\frac58.$

The blowing-up $B = Bl_Y \mathbb{A}^2$ is given by
$$B = \mbox{Proj} \frac{\C[x,y][A,B,C]}{(y^2A + x^3B, xC+yB,
B^3x+AC^2)}$$ and the exceptional variety $E$ consists of two
projective lines $E'$ and $E''.$ The map $\fie : \tilde{X} \to B$
maps $E_1$ and $E_3$ one-to-one onto $E'$ and $E'',$ respectively.
It contracts $E_2$ to a point $a$ on $E''.$

We choose a point $x \in E'^\circ$ and a point $y \in E''^\circ \<
\{a\}.$ For these points the monodromy zeta functions are
\[ Z_{\I , x}(t) = 1-t^3, \qquad  Z_{\I , y}(t) = 1-t^8.
\]
So we conclude that the two poles of the topological zeta function
induce two eigenvalues of monodromy $e^{-\frac{4 \pi i}{3}}$ and
$e^{-\frac{5 \pi i}{4}}.$

\begin{remark}
Note that this example is a monomial ideal. As mentioned in the
introduction, the monodromy conjecture had already implicitly been
verified in this case.
\end{remark}

\vspace{\baselineskip}

\begin{figure}[h]  \caption{\textsl{Intersection diagram of the principalization
of $(x^3y, x^3-y^2).$}}\label{principalisatie voorbeeld 3}
\unitlength=1mm
\begin{center}
\begin{picture}(90,32)(10,8)
\linethickness{0.5mm}

\put(32,17){\line(0,1){20}} \put(27,32){\line(1,0){46}}
\put(50,17){\line(0,1){20}} \put(68,12){\line(0,1){25}}
\put(63,17){\line(1,0){25}} \put(83,0){\line(0,1){22}}

\put(21,13){$E_1(2,2)$} \put(39,13){$E_2(3,3)$}
\put(14,30){$E_3(6,5)$} \put(69,37){$E_4(7,6)$}
\put(90,16){$E_5(8,7)$} \put(77,23){$E_6(9,8)$}

\linethickness{0.15mm}

\put(76,5){\line(1,0){3}} \put(81,5){\line(1,0){4}}
\put(87,5){\line(1,0){3}}

\end{picture}
\end{center}
\end{figure}
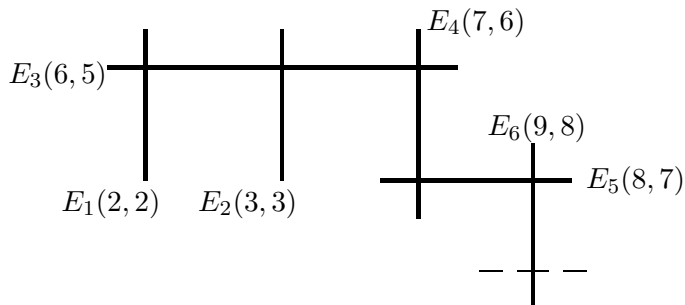

(3) $\I = (x^3y, x^3-y^2) \subset \C[x,y]$. In Figure
\ref{principalisatie voorbeeld 3} we can see that $Z_{top, \I}(s)$
has poles in $-\frac56$ and $-\frac89.$

We know that
$$B = Bl_YX = \mbox{Proj} \frac{\C[x,y][A,B]}{(A(x^3-y^2)-Bx^3y)}.$$
The blowing-up has one irreducible exceptional component $E.$ The
map $\fie : \tilde{X} \to B$ maps $E_6$ surjectively on $E$ and
all the other exceptional components of the principalization are
contracted to a point $a$ on $E.$ We fix a point $x \in E^\circ ,
x \neq a,$ and look at the following two monodromy zeta functions:
\[ Z_{\I , a}(t) = \frac{(1-t^2)(1-t^3)}{1-t^6}, \qquad  Z_{\I , x}(t) = 1-t^9.
\]
We can see that the two poles of the topological zeta function
$-\frac56$ and $-\frac89$ give rise to eigenvalues of monodromy
$e^{-\frac{5 \pi i}3}$ and $e^{-\frac{16 \pi i}9}.$

\footnotesize{
 
}

\end{document}